\documentclass[twocolumn]{autart}

\usepackage{amsmath}
\usepackage{amssymb}
\usepackage{makeidx}
\usepackage{multicol}
\usepackage[bottom]{footmisc}
\usepackage{graphicx}
\usepackage{color}
\usepackage{natbib}
\usepackage{epstopdf}
\usepackage{algorithm}
\usepackage{algcompatible}
\usepackage{calc}
\usepackage{subcaption}

\newtheorem{example}{Example}
\newtheorem{approxi}{Approximation}
\newtheorem{assump}{Assumption}

\newcommand{\app}{\begin{approxi}}
\newcommand{\eapp}{\end{approxi}}
\newcommand{\ass}{\begin{assump}}
\newcommand{\eass}{\end{assump}}
\newcommand{\teo}{\begin{thm}}
\newcommand{\eteo}{\end{thm}}
\newcommand{\corr}{\begin{cor}}
\newcommand{\ecorr}{\end{cor}}
\newcommand{\pro}{\begin{prop}}
\newcommand{\epro}{\end{prop}}
\newcommand{\lemma}{\begin{lem}}
\newcommand{\elemma}{\end{lem}}
\newcommand{\pb}{\begin{prob}}
\newcommand{\epb}{\end{prob}}
\newcommand{\df}{\begin{defn}}
\newcommand{\edf}{\end{defn}}
\newcommand{\rema}{\begin{rem}}
\newcommand{\erema}{\end{rem}}
\newcommand{\ex}{\begin{example}}
\newcommand{\eex}{\end{example}}

\newcommand{\al}[1]{\begin{align} #1 \end{align}}
\newcommand{\nn}{\nonumber}

\newcommand{\Zc}{ \mathcal{Z}}

\newcommand{\Cs}{ \mathbb{C}}

\newcommand{\Es}{ \mathbb{E}}

\newcommand{\Ns}{ \mathbb{N}}

\newcommand{\Rs}{ \mathbb{R}}

\newcommand{\Zs}{ \mathbb{Z}}

\begin{document}

\begin{frontmatter}
\title{The Harmonic Analysis of Kernel Functions \thanksref{footnoteinfo}}

\thanks[footnoteinfo]{This work has been partially supported by the FIRB project ``Learning
meets time'' (RBFR12M3AC) funded by MIUR.}

\author[Padova]{Mattia Zorzi}\ead{zorzimat@dei.unipd.it},
\author[Padova]{Alessandro Chiuso}\ead{chiuso@dei.unipd.it}

\address[Padova]{Dipartimento di Ingegneria dell'Informazione, Universit\`a degli studi di
Padova, via Gradenigo 6/B, 35131 Padova, Italy}


\begin{abstract} 
Kernel-based methods have been recently introduced for linear system identification as an alternative to parametric prediction error methods. 
Adopting the Bayesian perspective, the impulse response is modeled as a non-stationary Gaussian process with zero mean and with a certain kernel (i.e. covariance) function. Choosing the kernel is one of the most challenging 
and important issues. In the present paper we  introduce the harmonic analysis of this non-stationary process, and argue that this is an important tool which  helps in  designing such  kernel. Furthermore, this analysis suggests also an effective way to approximate the kernel, which allows to reduce the computational burden of the identification procedure. 
\end{abstract}
\end{frontmatter}

\section{Introduction}
Building upon the theory of reproducing kernel Hilbert spaces and statistical learning,  kernel-based methods for linear system identification  have been recently introduced in the system identification literature,  see  \cite{PILLONETTO_DENICOLAO2010,PILLONETTO_2011_PREDICTION_ERROR,CHIUSO_PILLONETTO_SPARSE_2012,EST_TF_REVISITED_2012,KERNEL_METHODS_2014,BSL,BSL_CDC,ChiusoARC2016,RPEM}. 

These methods, framed in the context of Prediction Error Minimization, differ from classical parametric methods  \cite{LJUNG_SYS_ID_1999,SODERSTROM_STOICA_1988},  in that models are searched for in possibly infinite dimensional model classes, described by Reproducing Kernel Hilbert Spaces (RKHS).  Equivalently, in a Bayesian framework, models are described assigning as prior  a Gaussian distribution;  estimation is then performed following the prescription of Bayesian Statistics, combining the ``prior'' information with the data in the posteriors distribution. Choosing the covariance function of the prior distribution, or equivalently the Kernel defining the RKHS, is one of the most challenging and important issues. For instance the prior distribution could reflect the fact that the system is Bounded Input Bounded Output (BIBO) stable, its impulse response  possibly smooth and so on \citep{KERNEL_METHODS_2014}. 

Within this framework, the purpose of this paper is to discuss the properties of certain kernel choices from the point of view of Harmonic Analysis of stationary processes. The latter is well defined for stationary processes  \citep[Chapter 3]{LINDQUIST_PICCI}. In particular, it defines as Power Spectral Density (PSD) the function describing how the statistical power is distributed over the frequency domain. In this paper, we extend this analysis for a particular class of non-stationary  processes modeling impulse responses of marginally stable systems. Accordingly, we define as Generalized Power Spectral Density (GPSD) the function describing 
how the statistical power  is distributed over the decay rate-frequency domain. 

Under the assumption that the prior density is Gaussian,  the probability density function (PDF) of the prior is linked to the GPSD. 
The main difference  is that while the former is defined over an infinite dimensional space (the underlying RKHS ${\mathcal H}_K$),  the latter is defined over the bidimensional  decay rate-frequency space. As a consequence, the latter is simple to depict but also to interpret from an engineering point of view. We show experimentally that, over the class of second-order linear systems, the two provide similar information. This class is important because: 1) it contains the simplest systems that exhibit oscillations and overshoot;  2) second order systems are building block of higher order systems and, as such,  understanding second order systems helps understanding higher ones. Furthermore, for a special class of exponentially convex locally stationary processes (ECLS) typically used in system identification \citep{KERNEL_STRUCTURES_II_2015,KERNEL_STRUCTURES_I_2015}, it is possible to provide (i) a clear link between the GPSD and the  Fourier transform of the exponentially modulated sample trajectory and (ii) characterize the posterior mean in terms of the GPSD. As a consequence,
it is possible to outline a simple procedure for the design of the kernel, through the GPSD.

Another important aspect in kernel-based methods is to reduce the computational burden \citep{Chen20132213}. This task can be accomplished by approximating the kernel functions through eigen-decomposition \citep{CaChPi2012} or random features, \citep{rahimi2007random} techniques. However, these methods can be applied only to special kernel functions.
We show that, the GPSD provides a general procedure to approximate a wide class  kernel functions.

The outline of the paper follows. In Section \ref{sec:gpr} we review the Gaussian process regression framework used for kernel-based methods. In Section \ref{sec_Harmonic_cont} we present the harmonic representation of the kernel function of continuos time non-stationary processes modeling impulse responses of marginally stable systems and in Section \ref{sec:Har_an_discrete} the corresponding discrete time version. In Section \ref{pdf_gpsd} we show the relation between the GPSD and the probability density function of the prior over the class of second-order linear systems.  In Section \ref{psd_ecls} we characterize the posterior mean in terms of the GPSD for a special class of ECLS processes.  Section \ref{sec:approx} regards kernel approximation using the GPSD.  Finally, conclusions are drawn in Section \ref{sec:concl}. In order to streamline the presentation all proofs are deferred to the Appendix.

{\em Notation}: $\Ns$ denotes the set of the natural numbers, $\Zs$ is the set of the integer numbers, $\Rs_+$ the set of the nonnegative real numbers, $\Rs_-$ the set of the negative real numbers, and $\Cs$ is the set of the complex numbers. Given $x\in\Cs$, $|x|$ denotes its absolute value, $\angle{x}$ denotes its phase and $\overline{x}$ denotes its conjugate. Given $A\in \Cs^{m\times n}$, $A^*$ denotes its transposed conjugate. $\Es[\cdot]$ denotes the expectation operator.

\section{System Identification and Gaussian Process Regression}\label{sec:gpr}

For convenience in what follows we consider a discrete time, stable and  linear time-invariant (LTI) single input-single output (SISO) in OE form: 
\al{\label{model}y(t)=G(z)u(t)+e(t),\;\; t\in\Ns}
where $z^{-1}$ denotes the backward shift operator; $u(t)$ is the input; $y(t)$ is the noisy output; $e(t)$ is zero mean white noise, that is $\Es[e(t)e(s)]=\sigma^2 \delta_{t-s}$  
where $\delta_k$ denotes the Kronecker delta function. The transfer function  $G(z)$ is  stable and strictly causal, i.e.  $G(\infty)=0$. 
Expanding $G(z)$ in $z^{-1}$ we obtain the impulse response of the linear system\al{G(z)&=\sum_{t=1}^\infty g(t) z^{-t}.\nn}  The system identification problem can be frased as that of  estimating the impulse response  $\{g(t)\}_{t\in\Ns}$, from the given data record 
 \al{ \Zc^N=\{ u(1),y(1) \ldots u(N),y(N)\}.\nn}
In the Gaussian process regression framework, $g(t)$ is modeled as a zero-mean discrete time Gaussian process with kernel (covariance) function $K(t,s): = \Es [g(t)g(s)]$. The minimum variance estimator of $g(t)$
is given by its posterior mean given $\Zc^N$,  \citep{PILLONETTO_DENICOLAO2010}. It is clear that the posterior  highly depends on the kernel functions. Accordingly, the most challenging part of this system identification procedure is to design  $K$ so that the posterior has some desired properties. 

Similarly, in the continuous time case,  $\{g(t)\}_{t\in \Rs_+}$ is  a zero-mean continuous time Gaussian  process with kernel function $K(t,s): = \Es [g(t)g(s)]$, with $t,s\in\Rs_+$. In what follows Gaussian processes (both discrete time and continuous time) are always understood with zero-mean.

\section{Harmonic Analysis: continuous time case}\label{sec_Harmonic_cont}
It is well known that the impulse response  of a finite dimensional LTI stable (or marginally stable) system 
can be written as a linear combination of decaying sinusoids\footnote{For simplicity we exclude here the case of multiple eigenvalues.} (i.e. modes)
\al{\label{def_g_t}g_a(t)=\sum_{l=1}^{N} |c_{l}|  e^{\alpha_l t}\cos(\omega_l t+\angle{c_{l}}),\;\; t\in\Rs_+}
where $\alpha_l \in\Rs_- \cup\{0\}$ and $\omega_l\in \Rs$ are, respectively, the decay rate and the angular frequency of the $l$-th damped oscillation, and  $c_l\in \Cs$. Adopting the Bayesian perspective, $g(t)$ is modeled as a Gaussian process  where the coefficients $c_l$  are zero mean complex Gaussian random variables such that
\al{& \Es[c_{l}\overline{ c_{l^\prime}}]=\phi_l \delta_{ l-l^\prime}\nn\\
& \Es[c_{l}c_{l^\prime}]=0\nn}
with $l,l^\prime=1\ldots N$ and $\phi_l\geq 0$. For convenience, we rewrite (\ref{def_g_t}) as 
\al{ \label{def_g_t2}g_a(t)=\sum_{i=1}^{N_\alpha}\sum_{k=1}^{N_\omega}  |c_{ik}|  e^{\alpha_{i}t }\cos(\omega_k t+\angle{c_{ik}})}
that is $(\alpha_i,\omega_k)$ belongs to a $N_\alpha \times N_\omega$ grid contained in $\Rs_- \cup\{0\}\times \Rs$ and 
$c_{ik}$ is a complex Gaussian random variable such that 
\al{& \Es[c_{ik}\overline{ c_{i^\prime k^\prime}}]=\phi_{ik}  \delta_{ i-i^\prime}\delta_{k-k^\prime}\nn\\
& \Es[c_{ik}c_{i^\prime k^\prime}]=0\nn} 
with $\phi_{ik}\geq 0$. It is then natural to generalize (\ref{def_g_t2}) as an infinite ``dense'' sum of decaying sinusoids \footnote{Strictly speaking the integral over $\alpha$ should be understood in an open interval of the form  $(-\infty, \epsilon)$ with $\epsilon>0$ and $\epsilon \rightarrow 0$.}:
\al{\label{def_g_t3}g(t)=\int_{-\infty}^0 \int_{-\infty}^{\infty} |c(\alpha,\omega )| e^{\alpha t} \cos(\omega t+\angle{c(\alpha,\omega)}) \mathrm{d}\omega \mathrm{d}\alpha }
 where $c(\alpha,\omega)$ is a bidimensional complex Gaussian process\footnote{To be precise, we should work with the ``generalized spectral measure'' $C(\omega,\alpha)$, a Gaussian  process with orthogonal increments; formally $dC(\omega,\alpha)=c(\omega,\alpha)d\omega d\alpha$.}, hereafter called generalized Fourier transform of $g(t)$, such that  
\al{& \Es[c(\alpha,\omega)\overline{ c(\alpha^\prime,\omega^\prime)}]=\phi(\alpha,\omega)  \delta (\alpha-\alpha^\prime) \delta (\omega-\omega^\prime)\nn\\
& \Es[c(\alpha,\omega)c(\alpha^\prime, \omega^\prime)]=0\nn}
where $\phi(\alpha,\omega)$ is a nonnegative function on $\Rs_-\cup\{0\}\times \Rs$ such that $\phi(\alpha,\omega)=\phi(\alpha,-\omega)$  and $\delta(\cdot)$  denotes the Dirac delta function.  
\pro \label{prop_har_repre}Let $K(t,s)$ be the kernel function of $g(t)$ in (\ref{def_g_t3}) then,
\al{\label{form_harm}K(t,s)=\frac{1}{2}\int_{-\infty}^0 \int_{-\infty}^{\infty} \phi(\alpha,\omega) e^{\alpha(t+s)}\cos(\omega (t-s))  \mathrm{d}\omega \mathrm{d}\alpha .}
\epro
Formula (\ref{form_harm}) is the harmonic representation of the covariance function of the non-stationary process (\ref{def_g_t3}). We refer to $\phi(\alpha,\omega)$ as generalized power spectral density (GPSD) of $g(t)$. The latter describes how the ``statistical power'' of $g(t)$ (which depends on $t$) is distributed over the decay rate-angular frequency space $\Rs_-\cup\{0\}\times \Rs$ according to: 
\al{\Es\left[ g(t)^2\right] = K(t,t)=\frac{1}{2}\int_{-\infty}^0 \int_{-\infty}^{\infty} \phi(\alpha,\omega)  e^{2\alpha t}\mathrm{d}\alpha  \mathrm{d} \omega. \nn }
In particular, when $t=0$, the exponential term $e^{2\alpha t}$ disappears, so that the statistical power of $g(0)$ is obtained as in the stationary case: 
\al{\label{K00}\Es\left[ g(0)^2\right] = K(t,t)=\frac{1}{2}\int_{-\infty}^0 \int_{-\infty}^{\infty} \phi(\alpha,\omega)  \mathrm{d}\alpha  \mathrm{d} \omega. }

It is worth noting that the GPSD can be understood as a function in $\Cs$, that is 
$\phi(s)$ with $s=\alpha+j\omega$, and its support is the left half-plane.
Now we show that the harmonic representation (\ref{form_harm}) describes many kernel functions used in system identification.

\subsection{Stationary kernels}
The special case of stationary processes  is recaptured when  $\phi(\alpha,\omega)=\delta(\alpha) \phi_1(\omega)$, with $\phi_1(\omega)=\phi_1(-\omega)$ nonnegative function; in fact, under this assumption, we have
\al{\label{ref_harm_st} K(t,s)=K_1(t-s)=\frac{1}{2}\int _{-\infty}^{\infty} \phi_1(\omega) \cos(\omega(t-s)) \mathrm{d} \omega}
which is a stationary kernel and $\phi_1(\omega)$ is the corresponding power spectral density (PSD). 
Note that, (\ref{ref_harm_st}) is the usual harmonic representation of a stationary  covariance function which is also exploited in spectral estimation problems \citep{TAC15,DUAL}.  
In this case the stationary process $g(t)$ is an infinite ``dense'' sum of sinusoids
\al{g(t)=\int_{-\infty}^{\infty} |c(\omega)| \cos(\omega t+\angle{c(\omega)}) \mathrm{d}\omega\nn}
which, more rigorously, should be written in terms of the spectral measure $C(\omega)$  \al{g(t)=\int_{-\infty}^\infty e^{j\omega t} \mathrm{d} C(\omega). \nn}
Since $\phi_1(\omega)$ is an even function, we can rewrite it as 
\al{\phi_1(\omega)=\tilde \phi_1(\omega)+\tilde \phi_1(-\omega)\nn}
with $\tilde \phi_1(\omega)$ nonnegative function. 
For instance, choosing $\tilde \phi_1(\omega)$ as one of the following: 
\al{\label{defpsdstaz}\tilde \phi_{L}(\omega)&= \frac{\beta}{\pi[\beta^2+(\omega-\omega_0)^2]} \nn\\
\tilde \phi_{C}(\omega)&= \frac{1 }{2\beta} e^{-\frac{|\omega-\omega_0|}{\beta}}\nn\\
\tilde \phi_{G}(\omega)&=\frac{1}{\sqrt{2\pi\beta }}  e^{-\frac{ (\omega-\omega_0)^2}{2\beta} } }
we obtain 
\al{K_{L}(t-s)&=\e^{-\beta |t-s|}\cos(\omega_0 (t-s))\nn\\
K_{C}(t-s)&=\frac{1}{1+\beta^2 (t-s)^2}\cos(\omega_0 (t-s))\nn\\
K_{G}(t-s)&=e^{-\frac{\beta}{2}\left(t-s\right)^2}\cos(\omega_0 (t-s))\nn}
where $\omega_0$ denotes the angular frequency for which $\tilde \phi$ takes the maximum and $\beta$ is proportional to the bandwidth. Setting $\omega_0=0$ we obtain, respectively, the Laplacian, Cauchy and Gaussian kernel \citep{RASMUSSEN_WILLIAMNS_2006}. In particular, the latter is widely used in robotics for the identification of the inverse dynamic \citep{ICUB_CDC}.

\subsection{Exponential Convex Locally Stationary (ECLS)  kernels}
A generalization of stationary kernels, introduced  by \cite{SILVERMAN_57}, is the so-called class of Exponentially Convex Locally Stationary (ECLS) kernels; this is obtained postulating a separable structure for the GPSD  $\phi(\alpha,\omega)=\phi_1(\omega)\phi_2(\alpha) $, with $\phi_1(\omega)$ and $\phi_2(\alpha)$ nonnegative functions, so that the kernel $K(t,s)$ inherits the decomposition 
\al{K(t,s)= K_1(t-s)K_2\left(t+s\right)\nn}
where $K_1$ has been defined in (\ref{ref_harm_st}) and 
\al{ K_2\left(t+s\right)=\int _{-\infty}^{0} \phi_2(\alpha) e^{\alpha(t+s)} \mathrm{d} \alpha.\nn}
 In the case that
$\phi_2\left( \alpha\right)=\delta(\alpha- \alpha_0)$ with $ \alpha_0\in \Rs_-$ we obtain the  ECLS kernel \citep{KERNEL_STRUCTURES_II_2015,chen2016kernel} :
\al{\label{ECLS_kernel}K(t,s)=e^{\alpha_0(t+s)} K_1(t-s).}
Specializing  $K_1(t-s)=K_L(t-s)$ with $\omega_0=0$, in (\ref{ECLS_kernel}) we obtain the so-called diagonal-correlated (DC) kernel. Furthermore,
for suitable choices of  $\phi_1(\omega)$ in (\ref{ECLS_kernel}), one can obtain the stable-spline (SS), the diagonal (DI) and the tuned-correlated kernel (TC), see \cite{EST_TF_REVISITED_2012}.

\subsection{Integrated kernels}
Consider the GPSD
\al{\phi(\alpha,\omega)=2\frac{-\alpha}{\pi[\omega^2+\alpha^2]  } \mathbf{1}_{[\alpha_m /2,\alpha_M/2]}(\alpha)\nn }
where \al{ \mathbf{1}_{[\alpha_m/2,\alpha_M/2]}(\alpha) = \left\{ \begin{array}{cc} 1, &  \alpha\in[\alpha_m /2,\alpha_M /2]  \\ 0, & \hbox{otherwise} \end{array}\right. \nn} with $\alpha_m<\alpha_M<0$. Then, the corresponding kernel function is 
\al{K(t,s)=&\int_{-\infty}^0   \mathbf{1}_{[\alpha_m /2,\alpha_M /2]}(\alpha) e^{\alpha (t+s)}\nn\\
& \times  \int_{-\infty}^{\infty} 
\frac{\label{intK_due}-\alpha}{\pi[\omega^2+\alpha^2]  }  \cos(\omega(t-s)) \mathrm{d} \omega \mathrm{d}\alpha\nn\\
& =\int_{\alpha_m/2}^{\alpha_M/2} e^{\alpha(t+s)} e^{\alpha|t-s|}\mathrm{d}\alpha\nn\\
& =\int_{\alpha_m/2}^{\alpha_M/2} e^{2\alpha\max\{t,s\}} \mathrm{d}\alpha\nn\\
&= \frac{e^{\alpha_M\max\{t,s\}}-e^{\alpha_m\max\{t,s\}}}{2\max\{t,s\} } 
}
which is similar to the integrated TC kernel introduced in \cite{Pillonetto2016137}  \footnote{See Remark \ref{oss_TCint} for more details.}. In general, taking 
\al{\phi(\alpha,\omega)=\phi_1(\omega;\alpha)  \mathbf{1}_{[\alpha_m,\alpha_M]}(\alpha),\nn  }
where we made explicit the dependence of $\phi_1$ upon $\alpha$, we have
\al{ \label{intK_uno} K(t,s)=\int_{\alpha_m}^{\alpha_M}e^{\alpha(t+s)} K_1(t-s;\alpha) \mathrm{d} \alpha  }
where $K_1(t-s;\alpha)$ is the stationary kernel corresponding to $\phi_1(\omega;\alpha)$. Notice that, 
kernel (\ref{intK_uno}) is obtained by integrating the ECLS kernel $e^{\alpha(t+s)} K_1(t-s;\alpha)$ over the interval $[\alpha_m,\alpha_M]$, which justifies the name ``integrated''. Another possible way to construct an integrated kernel is choosing 
\al{\phi(\alpha,\omega)=\phi_1(\omega) \mathbf{1}_{[\alpha_m,\alpha_M]}(\alpha)  \nn}
where $\phi_1$ does not depend on $\alpha$. Then, the corresponding kernel is 
\al{K(t,s)=K_1(t-s)\frac{e^{\alpha_M(t+s)}-e^{\alpha_m(t+s)}}{t+s}, \nn}
which is an ECLS kernel with 
\al{ K_2(t+s)=\frac{e^{\alpha_M(t+s)}-e^{\alpha_m(t+s)}}{t+s}. \nn}

\section{Harmonic analysis: discrete time case} \label{sec:Har_an_discrete}
Following the same argumentations of Section \ref{sec_Harmonic_cont}, a Gaussian process describing a discrete time 
causal impulse response can be understood as 
\al{g(t)=\int_0^1 \int_{-\pi}^\pi |c(\lambda,\vartheta)| \lambda^t \cos(\vartheta(t-s) ) \mathrm{d} \vartheta \mathrm{d}\lambda,\; \; t\in\Ns\nn}
where $c(\lambda,\vartheta)$ is the generalized Fourier transform of $g(t)$, $\lambda$ is the decay rate and $\vartheta$ is the normalized angular frequency. Moreover, the kernel function of $g(t)$ admits the following harmonic representation
\al{\label{harm_repr_discreto}K(t,s)=\frac{1}{2}\int_0^1 \int_{-\pi}^\pi \phi(\lambda,\vartheta) \lambda^{t+s} \cos(\vartheta(t-s)) \mathrm{d} \vartheta \mathrm{d}\lambda} and $\phi(\lambda,\vartheta)=\phi(\lambda,-\vartheta)$ is the GPSD of $g(t)$. The latter is a nonnegative function over the decay rate-normalized angular frequency space $[0,1]\times [-\pi,\pi]$. Also in this case the GPSD can be understood as a function in $\Cs$, that is $\phi(z)$ with $z=\lambda e^{j\vartheta}$, and its support is the unit circle.

In system identification, it is usual to design a kernel function for a continuous time Gaussian process $g_c(t)$, $t\in\Rs_+$.
Then, the ``discrete time'' kernel is obtained by sampling the ``continuous time'' kernel with a certain sampling time $T$. The latter corresponds to the discrete time Gaussian process $g_d(k)$, $k\in\Ns$, obtained by sampling $g_c(t)$. 
\pro \label{prop_sampledK}Let $\phi_c(\alpha,\omega)$ and $\phi_d(\lambda,\vartheta)$ be GPSD of $g_c(t)$, $t\in\Rs_+$, and $g_d(k)$, $k\in\Ns$, respectively.
Then,
\al{\label{eq_rep_per}\phi_d(\lambda,\vartheta)=\frac{1}{\lambda T^2}\sum_{k\in\Zs} \phi_c( T^{-1}\log\lambda, T^{-1}(\vartheta-2\pi k)).}
\epro
Accordingly, if the continuous time GPSD is such that  
\al{\label{cond_aliasing}\phi_c(\alpha,\omega)\approx 0,\;\;\; |\omega|>\frac{\pi}{T}}
then its discretized version is such that 
\al{\phi_d(\lambda,\vartheta)\approx \frac{1}{\lambda T^2}\phi_c( T^{-1}\log\lambda, T^{-1}\vartheta).\nn}
\rema \label{oss_TCint} Discretizing (\ref{intK_due}) with $T=1$, we obtain
\al{K_{iTC}(t,s)=\frac{\lambda_M^{\max\{t,s\}}-\lambda_m^{\max\{t,s\}}}{2\max\{t,s\}}, \nn}
with $0<\lambda_m<\lambda_M<1$ and $t,s\in\Ns$. However, the integrated TC kernel derived in \cite{Pillonetto2016137} is slightly different: 
\al{\label{KiTCbar}\bar K_{iTC}(t,s)=\frac{\lambda_M^{\max\{t,s\}+1}-\lambda_m^{\max\{t,s\}+1}}{\max\{t,s\}+1}.} Indeed, the latter has been derived by discretizing the TC kernel, and then the integration along the decay rate has been performed in the discrete domain. Note that,  the integration along the decay rate in $\bar K_{iTC}$is uniform, while, in $K_{iTC}$ such integration
is warped according to (\ref{eq_rep_per}).
\erema
\subsection{Filtered kernels}
In order to account for high frequency components of predictor impulse responses,  \cite{PILLONETTO_2011_PREDICTION_ERROR} have introduced a class of  priors obtained as filtered versions of  stable spline kernels, using second order filters of the form: 
\al{F(z)=\frac{z^2}{(z-\rho_0 e^{j\vartheta_0})(z-\rho_0 e^{-j\vartheta_0})} \nn}
with $|\rho_0|<1$ and $\vartheta_0\in[-\pi,\pi]$. The latter filter is fed by a Gaussian process $\tilde g(t)$ with kernel function $\tilde K(t,s)$, for instance in the class of ``stable-spline'' kernels (see \cite{pillonetto2010regularized}), which have most of the statistical power \eqref{K00} concentrated around  $\vartheta=0$; in this paper $\tilde K$ is chosen as TC kernel.  $F(z)$ plays the role of a shaping filter, which concentrates the  statistical power of the stationary part around $\vartheta=\vartheta_0$. It is not difficult to see that the kernel function of $g(t)$ is 
\al{K(t,s)=\sum_{l=1}^t \sum_{m=1}^s f_{t-l} f_{s-m} \tilde K(l,m) \nn}
where $\{f_l\}_{l\in \Ns}$ is the impulse response of the filter $F(z)$, i.e.  $F(z)=\sum_{s=0}^\infty f_sz^{-s}$. Assuming that $\tilde K(t,s)$ admits the harmonic representation 
\al{\tilde K(t,s)=\frac{1}{2}\int_0^1 \int_{-\pi}^\pi \phi(\lambda,\vartheta) \lambda^{t+s} \cos(\vartheta(t-s)) \mathrm{d} \vartheta \mathrm{d}\lambda\nn} then we have
\al{\label{harm_filtered}K(t,s)=\frac{1}{2}\int_0^1 \int_{-\pi}^\pi \phi(\lambda,\vartheta) f(\lambda,\vartheta,t,s) \mathrm{d} \vartheta \mathrm{d}\lambda}
where
\al{\label{filtered_basis}f(\lambda,\vartheta,t,s)=\sum_{l=1}^t \sum_{m=1}^s f_{t-l}f_{s-m}\lambda^{l+m} \cos(\vartheta(l-m)).}
Accordingly, (\ref{harm_filtered}) is the representation of the kernel function of $g(t)$ according to the basis functions in (\ref{filtered_basis}). Accordingly, the effect of filtering is to change the term
$\lambda^{t+s} \cos(\vartheta(t-s))$ with (\ref{filtered_basis}).

\section{Probability Density Function of Gaussian Processes  and their GPSD} \label{pdf_gpsd}
Given a discrete time Gaussian process  $g(t)$ with kernel $K(t,s)$,  in this Section we shall study the relation between the associated Gaussian probability distribution and the Generalized Spectral Density introduced in \eqref{form_harm}.  For simplicity we consider the discrete time case with sampling time $T=1$. In what follows, this process 
is understood as an infinite dimensional vector $g=[\,g(0)\,g(1)\,\ldots\,]^T$ and the corresponding kernel as an infinite dimensional matrix $K$. Then, the PDF of the process is 
\al{\label{pdf_g}p(g)\propto \exp\left( -\frac{g^T K^{-1} g}{2}\right)}
where ``$\propto$'' stands for ``proportional to''. Practically one can evaluate the r.h.s. of (\ref{pdf_g})  with respect to the impulse response $g$, which, with some abuse of notation, represents a realization of the process, truncating it according to its practical length and taking the corresponding finite dimensional sub-matrix\footnote{It can be shown that this makes sense provided $g$ belongs to the RKHS ${\mathcal H}$ with kernel $K$.} in $K$.
In this Section we consider a particular class of candidate impulse responses,
\begin{equation}\label{ClassSOS}\begin{array}{c}
 g \in {\mathcal G}_p:=\left\{g_p(t), t\in \Ns,\;\;\; s.t.\;\;\; G_p(z)=\frac{z^2}{(z-p)(z-\bar p)}\right.\\
\left. p\in \Cs \quad  |p|<1 \quad 0<\angle{p} <\pi\right\} 
\end{array}
\end{equation}
where 
$$
G_p(z):= \sum_{t=1}^{\infty} g_p(t) z^{-t}.
$$
Notice that each model in this class is uniquely characterized by $p$, accordingly this class is isomorphic to the upper part of the open unit circle. 

Thus, we consider the PDF \eqref{pdf_g} conditionally to the event $g \in {\mathcal G}_p$, 

\al{
p(g|g\in {\mathcal G}_p) \propto \begin{cases} \exp\left( -\frac{g^T K^{-1} g}{2}\right)   & g \in {\mathcal G}_p \\ 0 & g \notin {\mathcal G}_p. \end{cases}\nn
}

As we shall see, there is a close connection between the GPSD of the process $g(t)$ and $p(g|g\in {\mathcal G}_p)$, suggesting that indeed GPSDs whose stationary part  concentrate  energy around specific frequency bands are well suited to describe second order systems with modes in the same band (i.e. the phase of $p$ in \eqref{ClassSOS}). Similarly, the same applies to the decay rate, which relates to $|p|$ in  \eqref{ClassSOS}. 

 \subsection{ECLS kernels}
\begin{figure}[htbp]
\begin{center}
\includegraphics[width=\columnwidth]{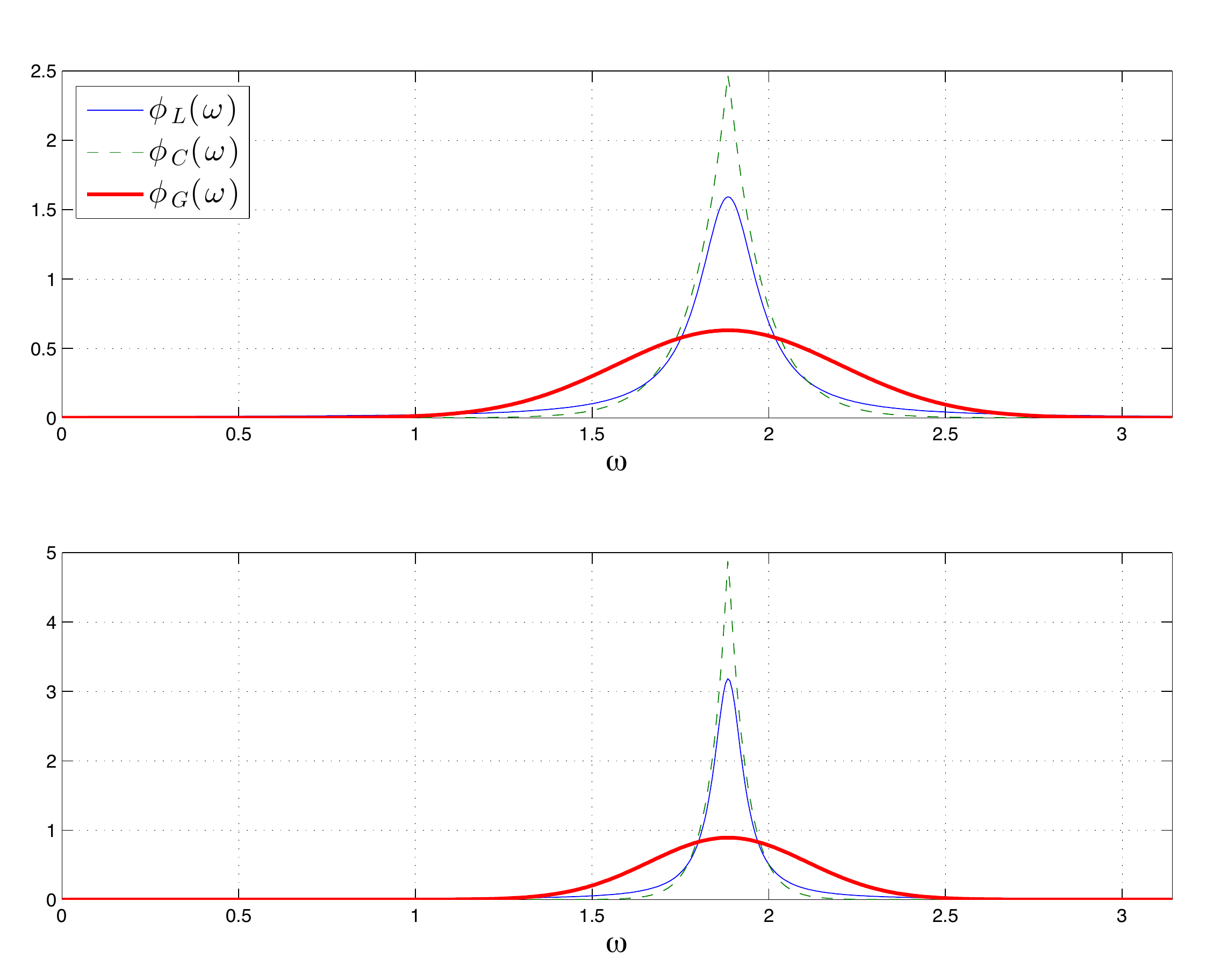}
\end{center}
\caption{PSD for the stationary part of the ECLS kernels  $K^{\hbox{\tiny ECLS}}_{L}$, $K^{\hbox{\tiny ECLS}}_{C}$ and $K^{\hbox{\tiny ECLS}}_{G}$ with  $\alpha_0=-0.1$, $\omega_0=3\pi/ 5$, $\beta=0.1$ (top) and $\beta=0.05$ (bottom).}\label{Fig1}
\end{figure}

\begin{figure}[htbp]
\begin{center}
\includegraphics[width=\columnwidth]{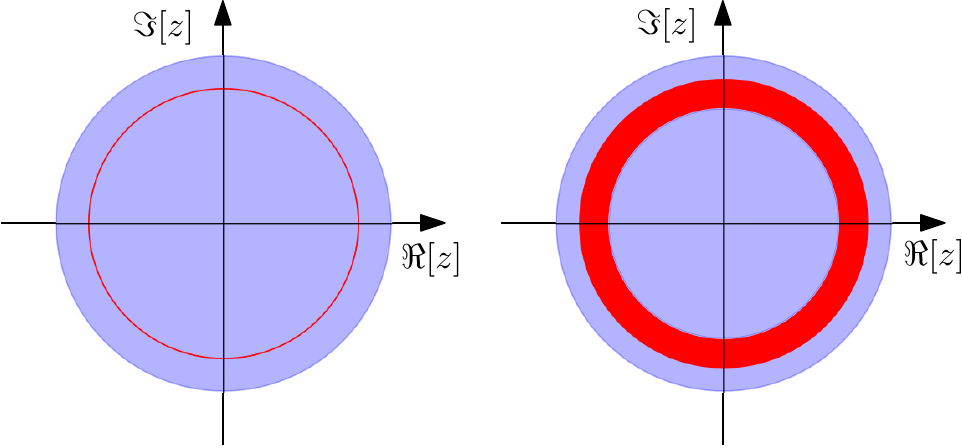}
\end{center}
\caption{Support of the GPSD $\phi(z)$ in red for an ECLS kernel (left) and for an integrated kernel (right).  The unit circle is depicted in transparent blue.}\label{Fig2}
\end{figure}

\begin{figure*}[htbp]
    \centering
            \begin{subfigure}{\textwidth}
            \centering
            \includegraphics[width=\textwidth]{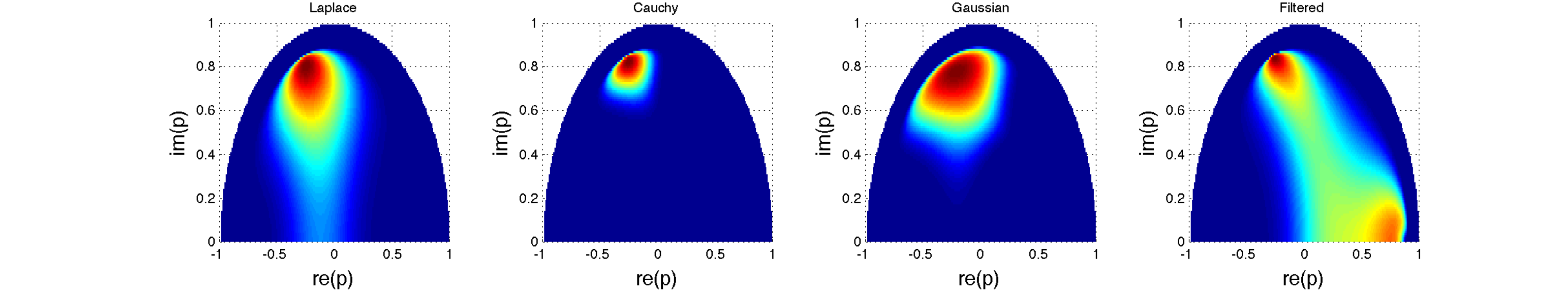}
            \caption{}
            \end{subfigure}\\
            \begin{subfigure}{\textwidth}
            \centering
            \includegraphics[width=\textwidth]{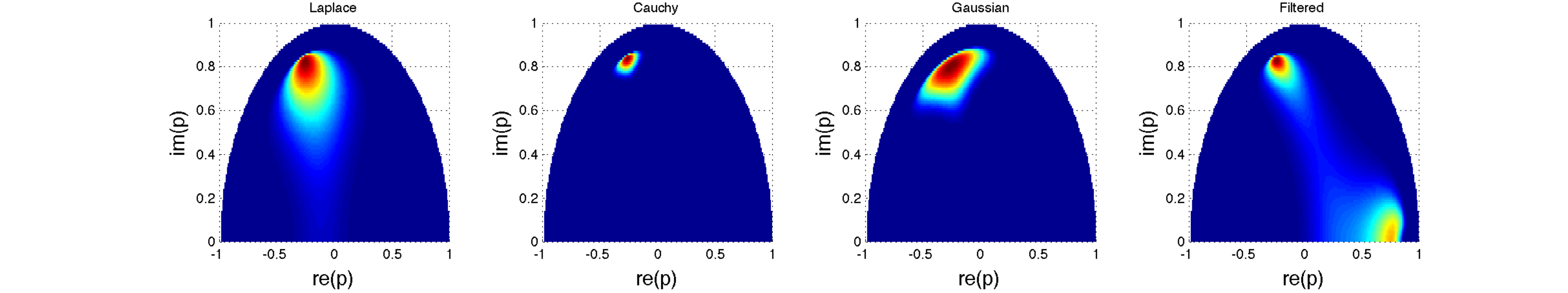}
            \caption{}
            \end{subfigure}%
            \\
            \begin{subfigure}{\textwidth}
            \centering
            \includegraphics[width=1.01\textwidth]{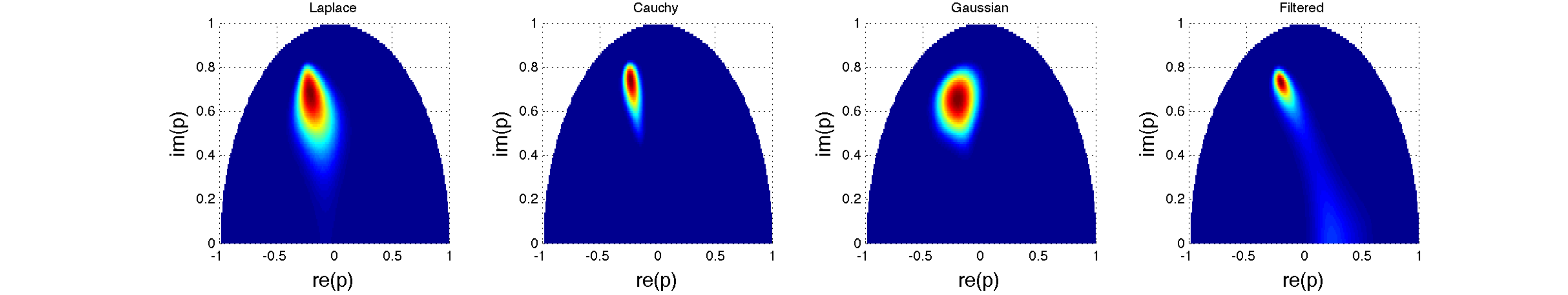}
            \caption{}
            \end{subfigure}%
            \caption{(a) PDF for processes having kernel $K^{\hbox{\tiny ECLS}}_{L}$ (first), $K^{\hbox{\tiny ECLS}}_{C}$ (second) and $K^{\hbox{\tiny ECLS}}_{G}$ (third) and filtered kernel (fourth). Here, $\beta=0.1$, $ \lambda_0=0.9$, $\vartheta_0=3\pi /5$ and $\rho_0=0.93$. (b) PDF for processes having kernel $K^{\hbox{\tiny ECLS}}_{L}$ (first), $K^{\hbox{\tiny ECLS}}_{C}$ (second) and $K^{\hbox{\tiny ECLS}}_{G}$ (third) and filtered kernel (fourth). Here, $\beta=0.05$, $ \lambda_0=0.9$, $\vartheta_0=3\pi /5$ and $\rho_0=0.96$. (c) PDF for processes having kernel $K^{\hbox{\tiny INT}}_{L}$ (first), $K^{\hbox{\tiny INT}}_{C}$ (second) and $K^{\hbox{\tiny INT}}_{G}$ (third) and integrated filtered kernel (fourth). Here, $\beta=0.05$, $\lambda_m=0.4$, $\lambda_M=0.9$, $\vartheta_0=3\pi /5$ and $\rho_0=0.96$.}\label{Fig3}
    \end{figure*}

We consider the ECLS kernels
\al{\label{ECLS_simulazione}K^{\hbox{\tiny ECLS}}_{L}(t,s)&=e^{ \alpha_0(t+s)}K_L(t-s)\nn\\
K^{\hbox{\tiny ECLS}}_{C}(t,s)&=e^{\alpha_0(t+s)}K_C(t-s)\nn\\
K^{\hbox{\tiny ECLS}}_{G}(t,s)&=e^{ \alpha_0(t+s)}K_G(t-s)}
where $\alpha_0=-0.1$, $\beta=0.1$ and $\omega_0=3\pi/ 5$. In Figure \ref{Fig1} (top) we show the corresponding PSDs of the stationary part.  
Then, we discretize these kernels with $T=1$. The corresponding GPSDs have as support a circle  with radius $\lambda_0=e^{\tilde \alpha_0 T}=0.9$ centered in zero, see left picture of Figure \ref{Fig2}.
Since condition (\ref{cond_aliasing}) holds,  the shape of the discretized versions  essentially reflect the continuous time version in  Figure \ref{Fig1} (top). 
 This means that the statistical power at time $t=0$ is concentrated in a neighborhood of $z_0=\lambda_0 e^{j\vartheta_0}$, with $\vartheta_0=\omega_0$, along the normalized angular frequency domain. In particular, the one of $K^{\hbox{\tiny ECLS}}_{C}$ is more concentrated in $z_0$ than the one of $K^{\hbox{\tiny ECLS}}_{L}$, and the latter is more concentrated than the one of  $K^{\hbox{\tiny ECLS}}_{G}$. Finally we also consider the filtered kernel with $\rho_0=0.93$ and $\vartheta_0=3\pi/ 5$. The corresponding PDFs
are depicted in Figure \ref{Fig3}(a).
  For all these kernels the impulse responses taking high probability have a pole in a neighborhood of $p_0= \lambda_0 e^{j\vartheta_0}$. Consistently with the GPSDs, the neighborhood of high-probability of $K_G^{ECLS}$ is more spread along $\vartheta$ than the one of $K_L^{ECLS}$, and the latter is more spread along $\vartheta$ than the one of $K_C^{ECLS}$. It is worth noting that for the filtered kernel there is a neighborhood of high-probability around a pole with phase close to zero. Indeed, process $\tilde g(t)$ at the input of $F(z)$ gives high-probability to impulse responses with pole whose phase is close to zero, and $F(z)$
does not sufficiently attenuate these impulse responses at the output.

Decreasing $\beta$ to $0.05$, the statistical power at time $t=0$ is more concentrated around $z_0$ along the normalized angular frequency domain, see Figure \ref{Fig1} (bottom). The corresponding PDFs are depicted in Figure \ref{Fig3}(b). Consistently with the GPSDs, the neighborhood of high-probability is more squeezed along $\vartheta$. Regarding the filtered kernel, $\rho$ has been increased to $0.96$; then the bandwidth of $F(z)$ is decreased. As a consequence, the neighborhood of high-probability around the pole with phase close to zero disappeared because $F(z)$ now drastically attenuate these impulse responses at the output.

\subsection{Integrated kernels}
We consider the integrated versions of (\ref{ECLS_simulazione})
\al{K^{\hbox{\tiny INT}}_{\mathcal{T}}(t,s)&= \frac{e^{\alpha_M(t+s)}-e^{\alpha_m(t+s)}}{t+s} K_{\mathcal{T}}(t-s)\nn}
where $\mathcal{T}\in\{ L,C,G\}$, $\alpha_m=-0.92$, $\alpha_M=-0.1$, $\beta=0.05$ and $\omega_0=3\pi/ 5$. Similarly as before, we discretize these kernels with $T=1$. The corresponding GPSDs have as support an annulus inside the unit circle, where the lower radius is $\lambda_m=e^{\alpha_m}=0.4$ and the upper radius $\lambda_M=e^{\alpha_M}=0.9$, see Figure \ref{Fig2} (right). 
If we develop any circle in the support, then we find the function of Figure \ref{Fig1} (bottom) up to a scaling factor. Therefore, the statistical power at time $t=0$ is concentrated in a neighborhood of $z_0$ which now is spread also along  the decay rate interval $[\lambda_m,\lambda_M]$. The corresponding PDFs are depicted in Figure \ref{Fig3} (c). Consistently with the GPSDs, the PDFs are such that  
 the neighborhood of high-probability probability around $p_0$ is now also spread over the decay rate domain. Note that, the integrated version of the filtered kernel is obtained by feeding $F(z)$, here $\rho_0=0.96$, with process $\tilde g(t)$ with kernel the integrated TC kernel (\ref{KiTCbar}).
 
\section{GPSD design for a special class of ECLS kernels} \label{psd_ecls}
In order to study in more detail the relation among the (sample) properties of Bayes estimators, the frequency domain description of the unknown impulse responses (Fourier transform)  and the properties of the kernel, we shall now specialize to a particular class of priors on $g(t)$, admitting an 
ECLS kernel as in (\ref{ECLS_kernel}). 
Under this restriction, 
the  generalized Fourier transform takes the form  
\al{\label{gFECLS}c(\alpha,\omega)=c_1(\omega)\delta (\alpha-\alpha_0)}
where $c_1(\omega)$ is a complex Gaussian process such that \al{\label{c_per_ecls} \Es[c_1(\omega)\overline{c_1(\omega^\prime)}]&=\phi_1(\omega)\delta (\omega-\omega^\prime)\nn\\   \Es[c_1(\omega){c_1(\omega^\prime)}]&=0,} and $\phi_1(\omega)\geq 0$ denotes its PSD. Indeed, 
\al{\Es[c(\alpha,\omega) &\overline{c(\alpha^\prime,\omega^\prime)}]\nn\\ & =\Es[c_1(\omega)\overline{c_1(\omega^\prime)}]\delta(\alpha-\alpha_0)\delta(\alpha^\prime-\alpha_0)\nn\\
&= \phi_1(\omega) \delta(\omega-\omega^\prime)\delta(\alpha-\alpha_0)\delta(\alpha-\alpha^\prime)\nn}
and $\phi(\alpha,\omega)=\phi_1(\omega)\delta(\alpha-\alpha_0)$.
By  (\ref{def_g_t3})  and (\ref{gFECLS}), it is not difficult to see that 
\al{\label{eq_gC} g(t)=\frac{1}{2}e^{\alpha_0 t} \int_{-\infty}^\infty G_{\alpha_0}(\omega)e^{j\omega t} \mathrm{d} \omega . }
where $G_{\alpha_0}(\omega):=(c_1(\omega)+\overline{c_1(-\omega)})/2 $.
Accordingly, $G_{\alpha_0}(\omega)$ is the Fourier transform of $\pi^{-1}e^{-\alpha_0 t} g(t)$:
\al{\label{trasf2}G_{\alpha_0}(\omega)=\int_{0}^\infty  \pi^{-1}e^{-\alpha_0 t} g(t) e^{-j\omega t} \mathrm{d} t.}
The next Proposition characterizes the posterior mean of $g(t)$ in terms of the PSD $\phi_1(\omega)$.   
\pro \label{prop_post_freq} Consider the continuos time model 
\al{\label{OEcont}v(t)&= \int_{0}^\infty g(s)u(t-s) \mathrm{d}s ,\; t\in\Rs_+}
 $u(t)$ is the (known) measured input. \\ Let $y^N=[\, y(1) \, \ldots \, y(N)\,]^T$ be sampled noisy measurements \footnote{For simplicity, here we assume that the sampling time is $T=1$.}
 $$
 y(k):= v(k) + e(k) \quad k=1,..,N
 $$
 where $e(k)$, $k=1,...,N$ are i.i.d. zero mean Gaussian with variance $\sigma^2$, independent of $g(t)$. 
 Assuming the prior distribution on $g(t)$ is Gaussian with kernel  (\ref{ECLS_kernel}), 
 the posterior mean  $\Es[g(t)|y^N]$ of $g(t)$ given $y^N=[\, y(1) \, \ldots \, y(N)\,]^T$ satisfies
 \al{\label{post_freq}\Es[g(t)|y^N]=\frac{1}{2}e^{\alpha_0 t}\int_{-\infty}^\infty \Es[G_{\alpha_0}(\omega)|y^N] e^{j\omega t}\mathrm{d}\omega} 
where \al{\label{post_freqsuG}\Es[G_{\alpha_0}(\omega)|y^N]=\phi_1(\omega) U_{\alpha_0}(\omega)^*V_{\alpha_0}^{-1} y^N} is the posterior mean of $G_{\alpha_0}(\omega)$ with 
\al{\label{def_Va0}V_{\alpha_0}&=\int_{-\infty}^\infty \phi_1(\omega)U_{\alpha_0} (\omega)U_{\alpha_0} (\omega)^* \mathrm{d}\omega  +\sigma^2 I_N\nn\\
U_{\alpha_0} (\omega)&=\frac{1}{2}\int_0^\infty e^{\alpha_0 s} u^N_s e^{j\omega s} \mathrm{d}s \nn\\
u^N_s &=[\, u(1-s) \, \ldots \, u(N-s)\,]^T.
} \epro
Proposition \ref{prop_post_freq} can be also adapted to the discrete time case as follows: 
\pro \label{prop_post_freq_discr} Consider a discrete time process $g(t)$ with ECLS kernel $K(t,s)=\lambda_0^{t+s}K_1(t-s)$, with $0<\lambda_0<1$ and $t,s\in\Ns$. Then, 
\al{g(t)&=\frac{1}{2}\lambda_0^t\int_{-\pi}^\pi G_{\lambda_0}(\vartheta)e^{j\vartheta t} \mathrm{d} \vartheta \nn\\ G_{\lambda_0}(\vartheta)&=\sum_{t=1}^\infty  \pi^{-1}\lambda_0^{-t} g(t) e^{-j\vartheta t} . \nn}
Consider the discrete time OE model 
\al{y(t)&=\sum_{s=1}^\infty g(s)u(t-s) + e(t)  ,\; t\in\Ns \nn}
where $u(t)$ is the measured input and $e(t)$ is zero mean white Gaussian noise with variance $\sigma^2$, independent of $g(t)$ and assume $\phi_1(\vartheta)$ is the PSD of $K_1(t-s)$.
Then the posterior mean   $\Es[g(t)|y^N]$  of $g(t)$ given $y^N=[\, y(1) \, \ldots \, y(N)\,]^T$ is given by
 \al{\Es[g(t)|y^N]=\frac{1}{2}\lambda_0 ^t\int_{-\pi
 }^\pi \Es[G_{\alpha_0}(\vartheta)|y^N] e^{j\vartheta t}\mathrm{d}\vartheta \nn} 
where \al{\Es[G_{\lambda_0}(\vartheta)|y^N]=\phi_1(\vartheta) U_{\lambda_0}(\vartheta)^*V_{\lambda_0}^{-1} y^N \nn} is the posterior mean of $G_{\lambda_0}(\vartheta)$ with 
\al{V_{\lambda_0}&=\int_{-\pi}^\pi \phi_1(\vartheta)U_{\lambda_0} (\vartheta)U_{\lambda_0} (\vartheta)^* \mathrm{d}\vartheta  +\sigma^2 I_N\nn\\
U_{\lambda_0} (\vartheta)&=\frac{1}{2}\sum_{s=1}^\infty e^{\lambda_0 s} u^N_s e^{j\vartheta s}  \nn\\
u^N_s &=[\, u(1-s) \, \ldots \, u(N-s)\,]^T.
\nn } \epro

For simplicity in what follows we consider the discrete time case, but the same observations hold for the continuous time case. Proposition \ref{prop_post_freq_discr}  shows that the absolute value of the  posterior mean of  $G_{\lambda_0}(\vartheta)$ is proportional to $\phi_1(\vartheta)$ (frequency wise). To understand better this fact,  we consider a data record $\mathcal{Z}^N$ of length $N=500$ generated by the discrete time model (\ref{model}) with transfer function having poles in $0.936$, $-0.45\pm 0.8$, $-0.25\pm 0.85$ and zeros in $0.16$, $-0.8\pm 0.4$. Since the dominant pole is $0.936$, we take $\lambda_0=0.94$.   \begin{figure}[htbp]
\begin{center}
\includegraphics[width=\columnwidth]{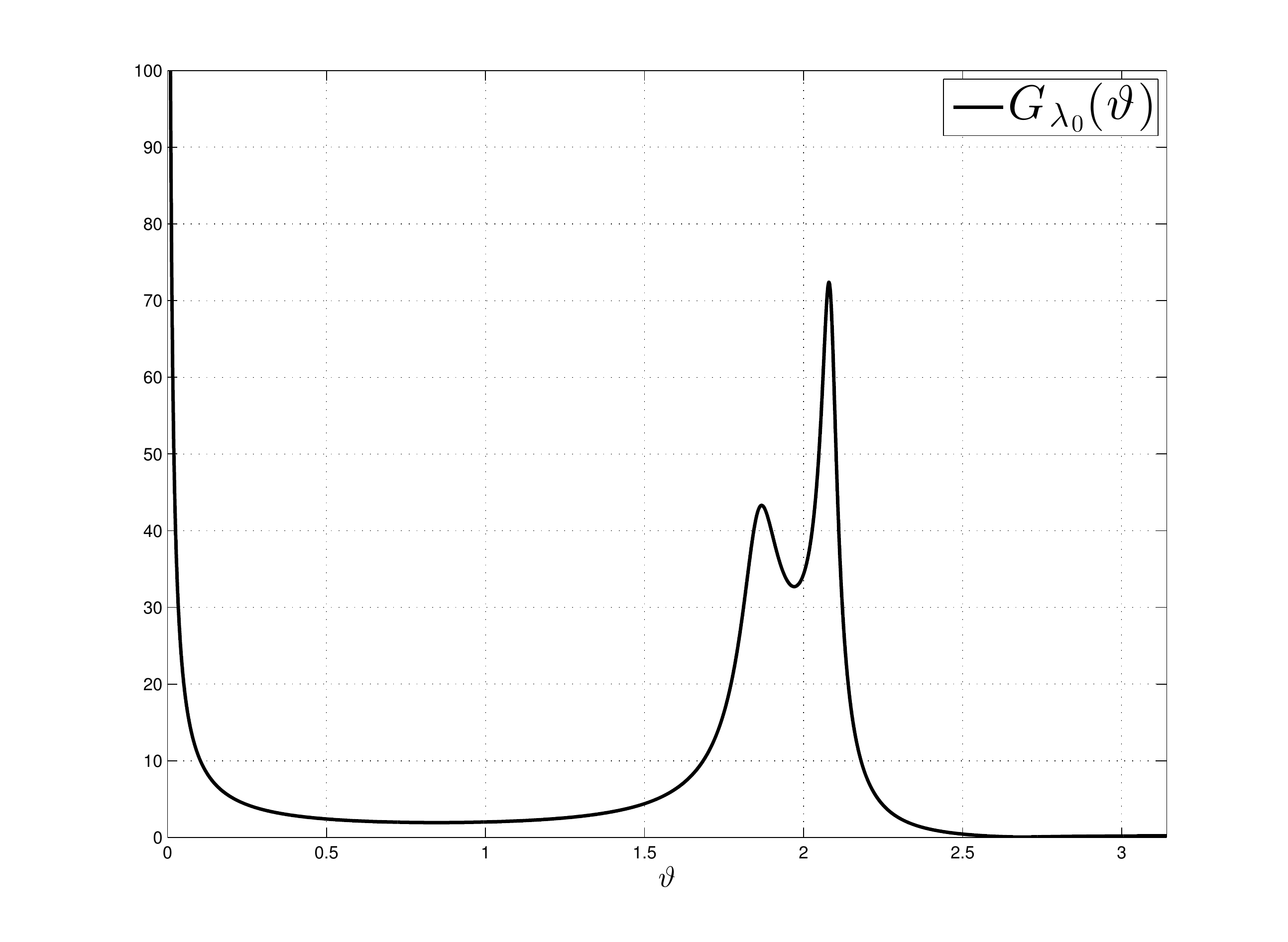}
\end{center}
\caption{Fourier transform of $\pi^{-1} \lambda_0^{-t} g(t)$.}\label{Figgtest}
\end{figure}
In this way $\pi^{-1}\lambda_0^{-t} g(t)$ admits Fourier transform $G_{\lambda_0}(\vartheta)$, see Figure \ref{Figgtest}.

We consider the  posterior mean estimators of $g(t)$ with three different sampled version of kernel (\ref{ECLS_kernel}) with $T=1$. In particular we choose $\alpha_0=\log (\lambda_0)$, and $\phi_1(\omega)$ takes one of the following three shapes:
\begin{itemize}
\item $\phi_1(\omega)=\phi_{L0}(\omega)=2\beta / (\pi(\beta^2+\omega^2))$
\item $\phi_1(\omega)=\phi_L(\omega)=\tilde\phi_L(\omega)+\tilde\phi_L(-\omega)$ with  $\tilde \phi_L$ as in (\ref{defpsdstaz})
\item $\phi_1(\omega)=\phi_M(\omega)=\phi_{L0}(\omega)+\phi_L(\omega)$ (a mixture kernel inspired by  \citep{MULTIPLE_KERNEL}).
\end{itemize}
We shall denote, correspondingly, with $g_{L0}(t)$, $g_{L}(t)$ and $g_{M}(t)$ the three estimators.

\begin{figure*}[htbp]
\begin{center}
\includegraphics[width=2\columnwidth]{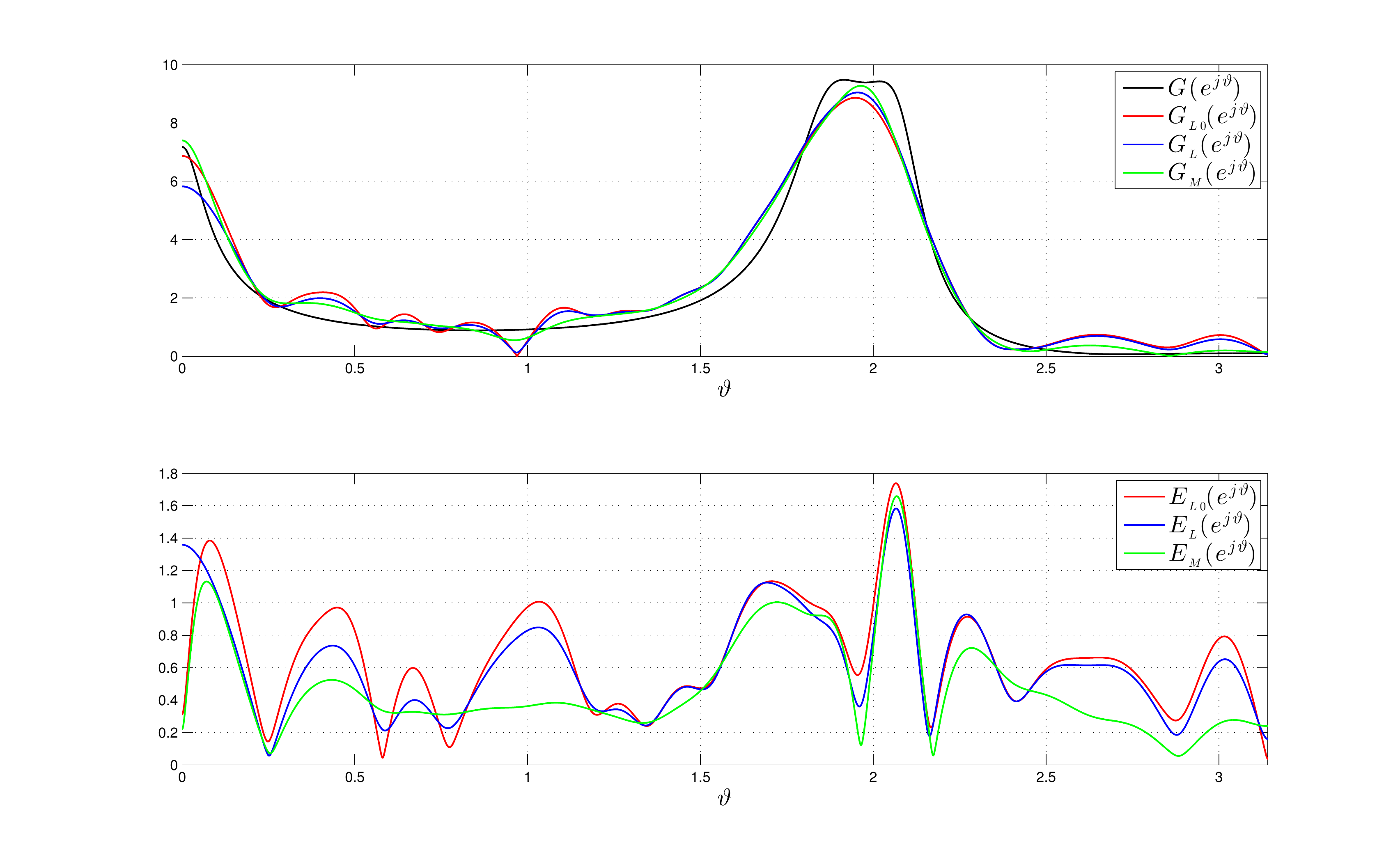}
\end{center}
\caption{Top panel: Fourier transform $G(e^{j\vartheta})$, $G_{L0}(e^{j\vartheta})$,  $G_L(e^{j\vartheta})$ and $G_M(e^{j\vartheta})$, respectively, of $g(t)$, $g_{L0}(t)$, $g_{L}(t)$ and $g_{M}(t)$. Bottom panel: absolute error of $G_{L0}(e^{j\vartheta})$,  $G_L(e^{j\vartheta})$ and $G_M(e^{j\vartheta})$ with respect to $G(e^{j\vartheta})$, e.g. $E_{L0}(e^{j\vartheta}):=|G_{L0}(e^{j\vartheta})-G(e^{j\vartheta})|$ .}\label{Figbode}
\end{figure*}
\begin{figure}[htbp]
\begin{center}
\includegraphics[width=\columnwidth]{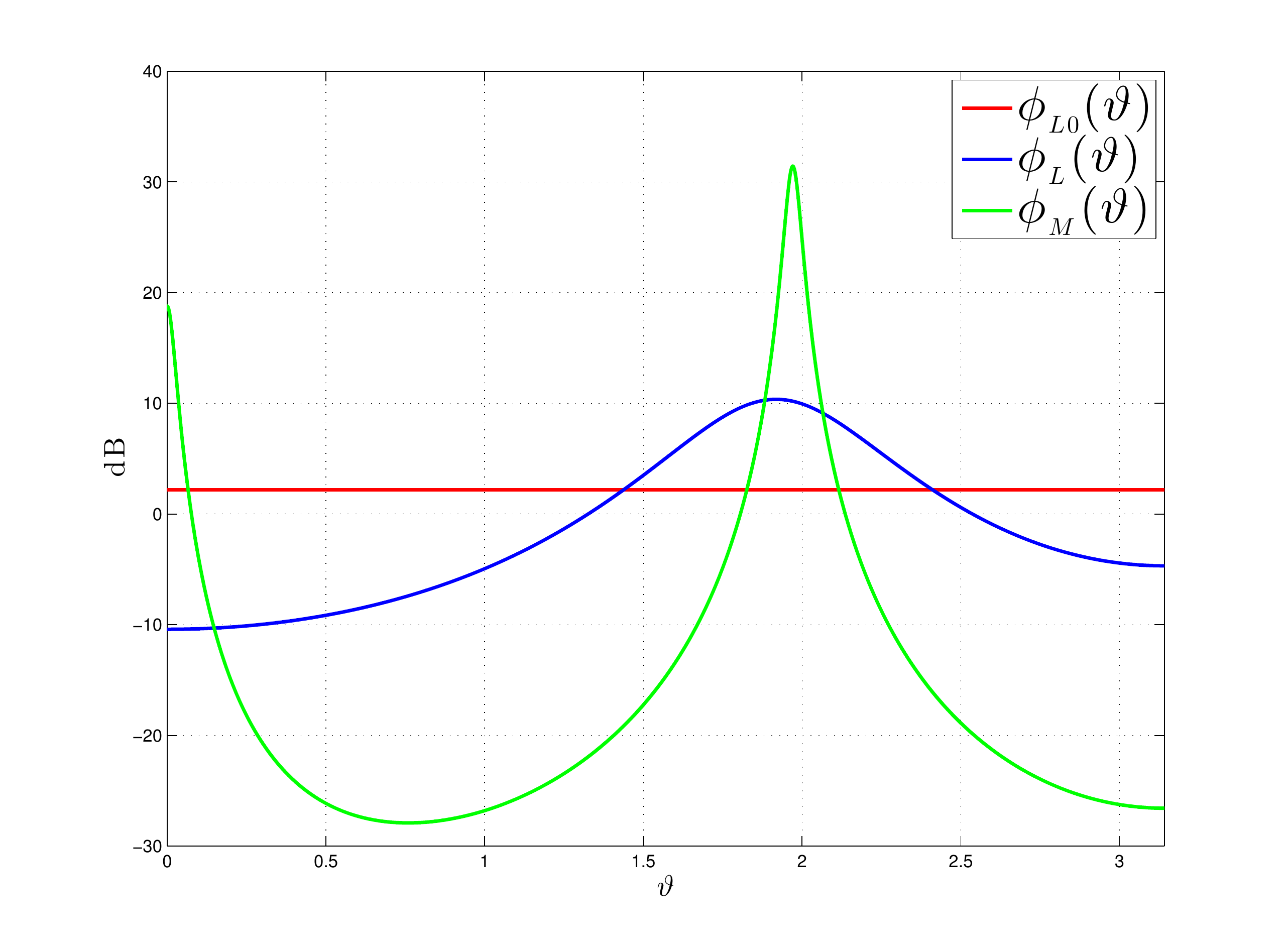}
\end{center}
\caption{PSD, in decibel, of the stationary part of the kernels corresponding to $g_{L0}(t)$, $g_L(t)$ and $g_{M}(t)$.}\label{Figgtest2}
\end{figure}

 The hyperparameters of the kernel are estimated from the data by minimizing the negative log-likelihood (e.g. the ones for $g_M(t)$ are $\beta$, $\omega_0$ and the scaling factor).
In Figure \ref{Figbode},
  the Fourier transform of the estimates, as well as the absolute errors are depicted. It is clear that $g_{M}(t)$ provides the best approximation of $g(t)$. In particular, compared to $g_{L0}(t)$ and $g_{L}(t)$, it improves the approximation 
at low ($\vartheta\simeq 0$), medium ($\vartheta\simeq 1$) and high ($\vartheta\simeq \pi$) normalized angular frequencies. In Figure \ref{Figgtest2},   the PSDs $\phi_{L0}(\vartheta)$, $\phi_{L}(\vartheta)$ and $\phi_{M}(\vartheta)$ of the stationary part of the three discretized kernels are depicted. Note that $\phi_{L0}(\vartheta)$, $\phi_{L}(\vartheta)$ and $\phi_{M}(\vartheta)$ are the periodic repetition (up to a scaling factor) of $\phi_{L0}(\omega)$, $\phi_{L}(\omega)$ and $\phi_{M}(\omega)$, respectively, with period $2\pi$.
It can be noticed that, only $\phi_{M}$ follows the shape of $|G_{\lambda_0}(\vartheta)|$.
This confirms the intuition that, if the PSD of the stationary part of the kernel has a similar shape of $|G_{\lambda_0}(\vartheta)|$, then we expect the corresponding estimate of $g(t)$ is good. 
 Hence this provides guidelines as to how the (stationary part of the) kernel should be designed; in particular, it should mimic the frequency response function $\mathrm{G}_{\lambda_0}^\dag(\vartheta)$ of  $\pi^{-1}\lambda_0^{- t} \mathrm{g}^\dag(t)$ where $\mathrm{g}^\dag(t)$ is the ``true'' impulse response function. 
 
 Note that  that optimality of the stationary part of the kernel is coupled with the choice of the decay $\lambda_0$; in practice, estimating  the hyperparameters using  the marginal likelihood estimator allows to optimize jointly $\lambda_0$ and the stationary part $\phi_1(\omega)$. 
 
 Note also that, for the discussion above to make sense, $\lambda_0^{- t} \mathrm{g}^\dag(t)$ should admit a Fourier transform, which imposes constraints on $\lambda_0$.
 To the purpose of illustration, let us postulate that  
  $\mathrm{g}^\dag(t)=\sum_{k=1}^\infty a_ie^{\alpha_i t}\cos(\delta_it)$ with $\alpha_i\in\Rs_-$, i.e. $\mathrm{g}^\dag(t)$ is a sum of damped sinusoids, then 
(\ref{eq_gC})-(\ref{trasf2}) hold if and only if \al{\label{conddec}\alpha_{max}:=\max_{i\geq 1} \alpha_i<\alpha_0:={\rm log}(\lambda_0).} 
For instance, the TC kernel, i.e. stable spline kernel of order one, is defined as $K_{TC}(t,s)=e^{-\gamma \max\{t,s\}}$ with $\gamma>0$ which can be written in the form  (\ref{ECLS_kernel})   choosing $\alpha_0=-\gamma /2 $ and $\phi_L(\omega)$ in (\ref{defpsdstaz}) with $\beta=\gamma /2$, therefore condition (\ref{conddec}) becomes $\alpha_{max}<-\gamma /2$. Furthermore, the SS kernel, i.e. stable spline kernel of order two, is defined 
\al{K_{SS}(t,s)=\frac{e^{-\gamma(t+s)} -e^{-\gamma \max\{t,s\}}}{2}-\frac{e^{-3\gamma\max\{t,s\}}}{6};\nn} we can rewrite it as (\ref{ECLS_kernel})  by choosing $\alpha_0=-3\gamma /2 $ and \al{\phi_1(\omega)=\frac{e^{-\frac{\gamma}{2}|t-s|}}{2}-\frac{e^{-\frac{3}{2}\gamma |t-s|}}{6},\nn} therefore condition (\ref{conddec})
becomes $\alpha_{max}< 3 \gamma /2$. It is interesting to note that the two conditions above coincide with the conditions derived in  \cite{pillonetto2010regularized} for  the posterior mean estimator to be statistically consistent, namely: 
 \al{\gamma< \frac{-2\alpha_{max}}{ 2m-1},\nn} where $m$ is the order of the stable spline kernel. Indeed, if (\ref{eq_gC})-(\ref{trasf2}) hold then 
  the hypothesis space is endowed by a probability density which is strictly positive at the ``true'' system  $g^\dag(t)$; under this condition it   can be proved that 
     the posterior mean estimator almost surely converges to  $\mathrm{g}^\dag(t)$. A similar  reasoning can be applied to the discrete time case.
%
%

\section{Kernel approximation} \label{sec:approx}
Kernel approximation is widely used in machine learning and system identification to reduce the computational burden. Next, we show that the GPSD represents a powerful tool for this problem.
To this purpose, note that the process $g_a(t)$ in (\ref{def_g_t2}) can be understood as an approximation of process $g(t)$ in (\ref{def_g_t3}). 
In particular, its kernel function is
\al{K_a(t,s)=\frac{1}{2}\sum_{i=1}^{N_\alpha} \sum_{k=1}^{N_\omega}  \phi_{ik} e^{\alpha_i(t+s)} \cos(\omega_k(t-s)).\nn } 
We define $z_{\alpha,\omega}(t)=e^{\alpha t}[\, \cos(\omega t)\;\, \sin(\omega t)]^\top$, then
\al{\label{def_Ka}K_a(t,s)=\frac{1}{2}\sum_{i=1}^{N_\alpha} \sum_{k=1}^{N_\omega}  \phi_{ik} z_{\alpha_i,\omega_k}(t)^\top z_{\alpha_i,\omega_k}(s).}
The computational burden of  the identification procedure based on $K_a(t,s)$ is related to the number of points $N_\alpha$ and $N_\omega$, which can thus be chosen to trade off kernel approximation and computational cost. For fixed $N_\alpha$ and $N_\omega$, the quality of the approximation   $K_a(t,s)\approx K(t,s)$ depends on the choice of  the  $\phi_{ik}$'s, which will be now discussed.  First of all, let us observe that (\ref{K00}) can be approximated as follows:
\al{\label{K00a}K(0,0)\approx \frac{1}{2} \sum_{i=1}^{N_\alpha} \sum_{k=1}^{N_\omega}  \Delta \alpha_i \Delta \omega_k \phi(\alpha_i,\omega_k)}
where $\Delta\alpha_i=\alpha_{i+1}-\alpha_i$, $\Delta\omega_k=\omega_{k+1}-\omega_k$ with $\alpha_{i+1}>\alpha_i$, $\omega_{k+1}>\omega_k$. On the other hand, we have
\al{\label{Ka00}K_a(0,0)=\frac{1}{2}\sum_{i=1}^{N_\alpha} \sum_{k=1}^{N_\omega} \phi_{ik}.}
Matching (\ref{K00a}) and (\ref{Ka00}), we obtain
\al{\phi_{ik}=\Delta \alpha_i \Delta \omega_k \phi(\alpha_i,\omega_k). \nn}

An alternative way to approximate $K(t,s)$ takes inspiration from the random features approach for stationary kernels \citep{rahimi2007random}. Observe that (\ref{form_harm}) can be rewritten as 
\al{\label{K00E}K(t,s)=\frac{K(0,0)}{2}\Es[z_{\alpha,\omega}(t)^\top z_{\alpha,\omega}(s)  ]}
where $\Es[\cdot]$ is the expectation operator taken with respect to the PDF
$\tilde p(\alpha,\omega)=K(0,0)^{-1} \phi(\alpha,\omega)$. Then, we can approximate (\ref{K00E}) with
\al{K_a(t,s)=\frac{\gamma}{2}\sum_{i=1}^{N_\alpha} \sum_{k=1}^{N_\omega}  \phi_{ik} z_{\alpha_i,\omega_k}(t)^\top z_{\alpha_i,\omega_k}(s) \nn} 
where $\phi_{ik}$'s are drawn from $\tilde p(\alpha,\omega)$ and the constant $\gamma$ satisfies $\gamma=\frac{1}{2}\sum_{i=1}^{N_\alpha} \sum_{k=1}^{N_\omega}  \phi_{ik}$.
 
\subsection{Numerical experiments} 
\textbf{Data set}. We generate 1000 discrete time SISO systems of  order $30$. The poles are randomly generated as follows:  75\%  have phase randomly generated over an interval of size $\pi / 6$ centered in $\vartheta_0 \sim {\mathcal U}[\pi/ 4,3/ 4 \pi ]$ and absolute value $\sim {\mathcal U}[0.8, 0.95]$; the remaining poles are generated uniformly inside the closed unit disc of radius $[0, 0.95]$. For each system a data set of  230  points is obtained feeding the linear system with zero mean, unit variance, white Gaussian noise and corrupting the output with   additive zero mean white Gaussian noise so as to guarantee a signal to noise ratio equal to $10$. 
 
\textbf{Simulation setup and results}. We consider model (\ref{model}) with $G(z)=\sum_{t=1}^n g(t) z^{-t}$ where  $n=100$ is the practical length. 
We consider several estimators of $g(t)$ which differ on the choice of the kernel describing the prior distribution on $g(t)$. In particular we shall use the following subscripts:
\begin{itemize}
\item  $L0$ for estimator which uses the kernel $K^{INT}_L$ with  $\omega_0=0$;
\item  $L0A$ as above but with  $K^{INT}_L$   approximated using   (\ref{def_Ka}) with $N_\omega=5$ and $N_\alpha=3$; 
\item $L$ for estimator which uses the kernel $K^{INT}_L$;
\item $LA$ as above but with  $K^{INT}_L$   approximated using   (\ref{def_Ka}) with $N_\omega=5$ and $N_\alpha=3$;
\end{itemize}
The subscripts  $G0$, $G0A$, $G$, $GA$ and $C0$, $C0A$, $C$, $CA$ will have the same meaning as above but with respect to kernel $K^{INT}_G$ and $K^{INT}_C$ respectively.
Note that all hyperparameters (e.g. $\alpha_m$, $\alpha_M$, $\beta$, $\omega_0$ and the scaling factor for $K^{INT}_L$)  are estimated from the data by minimizing the negative log-likelihood. To evaluate the various kernels, the impulse response estimates $\hat g(t)$, $t=1\ldots n$, are compared to the true one, i.e. $g(t)$, by the average fit
\begin{align*}
AF=100\left(1-\sqrt{\frac{\sum_{t=1}^n |\hat g(t)-g(t)|^2} {\sum_{t=1}^n |\hat g(t)-\bar g|^2}}\right), \;\; \bar g=\frac{1}{n}\sum_{t=1}^n  g(t).
\end{align*}
\begin{figure}[htbp]
\begin{center}
\includegraphics[width=\columnwidth]{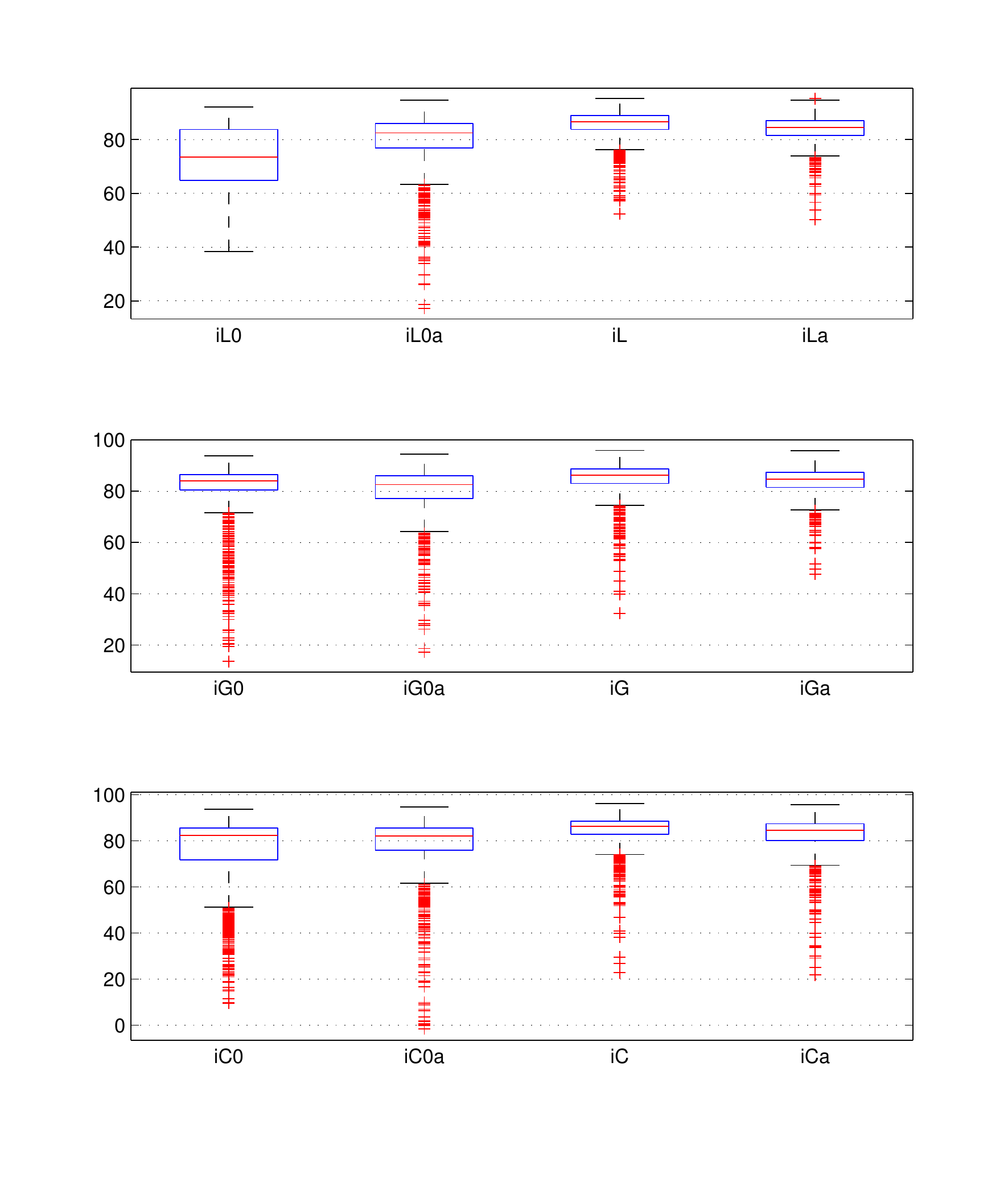}
\end{center}
\caption{Box-plots of the 1000 average fits.}\label{Fig6}
\end{figure} The distribution of the fits are shown by box-plots in Figure \ref{Fig6}: $L0$ and $C0$ estimators are outperformed by their approximated versions. In the remaining cases the approximated versions provides similar performance to the exact version. 

\section{Conclusions} \label{sec:concl}
In this paper we have introduced the harmonic representation of kernel functions used in system identification. In doing that, we have introduced the GPSD which represents the generalization of the power spectral density used for the harmonic representation of stationary kernels. We have showed by simulation that the GPSD carries similar information that we can find in the PDF of the process over the class of second-order systems. Moreover, we have characterized the posterior mean in terms of the GPSD for a special class of ECLS kernels. Finally, we have showed that the GPSD provides a powerful tool to approximate kernels function, and thus to reduce the computational burden in the system identification procedure.

\section*{Appendix} 

 \subsection*{Proof of Proposition \ref{prop_har_repre}}
 Rewriting 
\al{g(t)=\frac{1}{2}\int_{-\infty}^0\int_{-\infty}^{\infty}e^{\alpha t}  [c(\alpha,\omega) e^{j\omega t}+\overline{ c(\alpha,\omega)} e^{-j\omega t}] \mathrm{d} \omega \mathrm{d} \alpha\nn}  
then we have
\al{K&(t,s)=\Es[g(t)g(s)]\nn\\
&=\frac{1}{4}  \int_{-\infty}^{\infty}\int_{-\infty}^{\infty} \int_{-\infty}^0\int_{-\infty}^{0}e^{\alpha t+\alpha^\prime s}\nn\\ 
& \times \{ \Es[c(\alpha,\omega) c(\alpha^\prime,\omega^\prime)] e^{j\omega t+j\omega^\prime s}\nn\\
&+\Es[c(\alpha,\omega)\overline{ c(\alpha^\prime,\omega^\prime)}] e^{j\omega t-j\omega^\prime s}\nn\\
&+\Es[\overline{ c(\alpha,\omega)}c(\alpha^\prime,\omega^\prime)] e^{-j\omega t+ j\omega^\prime s}\nn\\
&+\Es[\overline{ c(\alpha,\omega)}\,\overline{ c(\alpha^\prime,\omega^\prime)}] e^{-j\omega t-j \omega^\prime s} \}\mathrm{d}\alpha\mathrm{d}\alpha^\prime \mathrm{d} \omega\mathrm{d}\omega^\prime \nn\\
&=\frac{1}{4}  \int_{-\infty}^{\infty}\int_{-\infty}^{\infty} \int_{-\infty}^0\int_{-\infty}^{0} e^{\alpha t+ \alpha^\prime s}   \phi(\alpha,\omega) \delta(\alpha-\alpha^\prime)\nn\\
&\times \delta(\omega-\omega^\prime) [ e^{j\omega t-j\omega^\prime s}+ e^{-j\omega t+ j\omega^\prime s} ] \mathrm{d}\alpha \mathrm{d}\alpha^\prime \mathrm{d} \omega\mathrm{d}\omega^\prime \nn\\
&=\frac{1}{4}  \int_{-\infty}^{\infty} \int_{-\infty}^{0}  e^{\alpha(t+s)}   \phi(\alpha,\omega) [ e^{j\omega (t-s)}+ e^{-j\omega (t-s)} ] \mathrm{d}\alpha  \mathrm{d} \omega\nn\\
&=\frac{1}{2}  \int_{-\infty}^{\infty}\int_{-\infty}^0  e^{\alpha(t+s)}   \phi(\alpha,\omega) \cos(\omega (t-s)) \mathrm{d}\alpha\mathrm{d} \omega. \nn  }
\qed

\subsection*{Proof of Proposition \ref{prop_sampledK}}
Let $K(t,s)$ be the kernel of $g_c(t)$. Define \al{\label{hts}h_{t+s}(t-s):=K(t,s).} Then, by (\ref{form_harm}) we have
\al{h_s(t)&=\frac{1}{2}\int_{-\infty}^0 \int_{-\infty}^{\infty} \phi_c(\alpha,\omega) e^{\alpha s}e^{j\omega t}  \mathrm{d}\omega \mathrm{d}\alpha \nn\\
&=\frac{1}{2\pi } \int_{-\infty}^{\infty}\underbrace{\int_{-\infty}^0 \pi \phi_c(\alpha,\omega) e^{\alpha s} \mathrm{d}\alpha}_{:= H_s(\omega)} e^{j\omega t}  \mathrm{d}\omega, \nn }
accordingly $ H_s(\omega)$ is the Fourier transform of $h_s(t)$. Now, we consider $h_s^d(k):=h_s(kT)$, $k\in\Ns$, which is the sampled 
version of $h_s(t)$ and $T$ is the sampling time. Therefore,
\al{\label{hk}h_s^d(k)=\frac{1}{2\pi}\int_{-\pi}^{\pi} H_s^d(\vartheta) e^{j\vartheta  k} \mathrm{d}\vartheta}
where $\vartheta=\omega T$ is the normalized angular frequency 
and 
\al{H_s^d(\vartheta):=\frac{1}{T}\sum_{l\in\Zs} H_s(T^{-1}(\vartheta-2\pi l )). \nn} 
Setting $n\in\Ns$, we have
\al{ \label{Hbar}&  H_{nT}^d (\vartheta) := \frac{1}{T}\sum_{l\in\Zs}H_{nT}(T^{-1}(\vartheta-2\pi l))\nn\\
&= \int_{-\infty}^0 \frac{\pi}{T} \sum_{l\in\Zs} \phi_c(\alpha,T^{-1}(\vartheta-2\pi l))e^{\alpha T n} \mathrm{d}\alpha\nn\\
&= \int_{0}^1 \pi \underbrace{\frac{1}{\lambda T^2}\sum_{l\in\Zs} \phi_c(T^{-1}\log \lambda,T^{-1}(\vartheta-2\pi l))}_{=\phi_d(\lambda,\vartheta)}\lambda^{ n} \mathrm{d}\lambda } where we used the substitution $\lambda=e^{\alpha T}$.
Substituting (\ref{Hbar}) in (\ref{hk}) and in view of (\ref{hts}), we obtain 
\al{K(kT,nT)=\frac{1}{2}\int_{0}^1 \int_{-\pi}^{\pi} \phi_d(\lambda,\vartheta) \lambda^{k+ n}e^{j\vartheta (k-n)}  \mathrm{d}\vartheta \mathrm{d}\lambda. \nn} In view of the discrete time harmonic representation 
(\ref{harm_repr_discreto}), we conclude that $\phi_d$ is the GPSD of the sampled kernel $K(kT,nT)$ corresponding to process $g_d$.\qed

\subsection*{Proof of Proposition \ref{prop_post_freq}}
By (\ref{eq_gC}), we trivially have (\ref{post_freq}). Then,
\al{\label{eq_part_prop2} \Es[G_{\alpha_0}(\omega)|y^N]&=\Es\left[\left.\frac{c_1(\omega)+\overline{c_1(-\omega)}}{2}\right|y^N\right]\nn \\&= \frac{\Es[c_1(\omega)|y^N]+\Es[\overline{c_1(-\omega)}|y^N]}{2}}
where
\al{\Es[c_1(\omega)|y^N]&=\Es[c_1(\omega)(y^N)^T]\Es[y^N (y^{N})^{ T}]^{-1}y^N\nn\\ 
\Es[\overline{c_1(-\omega)}|y^N]&=\Es[\overline{c_1(-\omega)}(y^N)^T] \Es[y^N (y^{N})^{ T}]^{-1}y^N.\nn} By (\ref{OEcont}), we have  
\al{v(t)
&=  \int_{-\infty}^\infty \frac{c_1(\omega)+\overline{c_1(-\omega)}}{2}U_{\alpha_0,t}(\omega) \mathrm{d}\omega, \nn} where we exploited relation (\ref{eq_gC}) and 
\al{U_{\alpha_0,t}(\omega):=\frac{1}{2}\int_0^\infty e^{\alpha_0 s} u(t-s)e^{j\omega s} \mathrm{d}s\nn.}  Let $k,n=1\ldots N$. It follows that 
\al{\Es[c_1(\omega)y(k)]&=\int_{-\infty}^\infty  \frac{\Es[c_1(\omega)c_1(\omega^\prime)]+\Es[{c_1(\omega)}\overline{c_1(-\omega^\prime)}]}{2} \nn\\ & \hspace{0.3cm} \times U_{\alpha_0, k}(\omega^\prime) \mathrm{d} \omega^\prime+\Es[c_1(\omega) e(k)]\nn\\ &=\frac{1}{2}\phi_1(\omega) \overline{U_{\alpha_0,k}(\omega) \nn}
} where we exploited (\ref{c_per_ecls}) and the fact that $U_{\alpha_0,t}(\omega)$ is an Hermitian function. Moreover,
\al{\label{eqEyts}\Es [y(k)y(n)]=\int_0^\infty\int_0^\infty & g(t^\prime) u(k-t^\prime)g(s^\prime) u(n-s^\prime) \mathrm{d}t^\prime\mathrm{d}s^\prime\nn\\& +\sigma^2 \delta_{k-n}.}
Substituting (\ref{eq_gC}) in (\ref{eqEyts}) and using (\ref{c_per_ecls}), we obtain
\al{\Es [y(k)y(n)]=\frac{1}{2}&\int_{-\infty}^\infty \phi_1(\omega^\prime) U_{\alpha_0,k}(\omega^\prime) \overline{U_{\alpha_0,n}(\omega^\prime)} \mathrm{d}\omega^\prime\nn\\ &+\sigma^2 \delta_{k-n}\nn} where we exploited the fact that $\phi_1(\omega)=\phi_1(-\omega)$. Since $U_{\alpha_0}(\omega)=[\, U_{\alpha_0,1} (\omega) \,\ldots\,\,U_{\alpha_0,N}(\omega)]^T$, we have that $\Es[c_1(\omega)(y^N)^T]=\phi_1(\omega) U_{\alpha_0}(\omega)^* / 2$ and
\al{\label{Ec1def} \Es[c_1(\omega)|y^N]=\phi_1(\omega) U_{\alpha_0}(\omega)^* V_{\alpha_0}^{-1} y^N}
where $V_{\alpha_0}$ has been defined in (\ref{def_Va0}). In similar way, it can be proven that \al{\label{Ec2def}\Es[\overline{c_1(-\omega)}|y^N]=\phi_1(\omega) U_{\alpha_0}(\omega)^* V_{\alpha_0}^{-1} y^N.} Finally, substituting (\ref{Ec1def}) and (\ref{Ec2def}) in (\ref{eq_part_prop2}) we obtain (\ref{post_freqsuG}).
\qed
\subsection*{Proof of Proposition \ref{prop_post_freq_discr}}
The proof is similar to the one of Proposition \ref{prop_post_freq}. \qed

\end{document}